\documentclass[11pt]{amsart}
\hoffset= - 0.75 in

\usepackage{amsfonts, amsmath, amssymb, amsthm}
\usepackage{latexsym, graphicx, color}
\newtheorem{thm}{Theorem}[section]

\newtheorem{defn}[thm]{Definition}

\newtheorem{prop}[thm]{Proposition}

\newcommand{\QQ}{\mathbb{Q}}
\newcommand{\RR}{\mathbb{R}}

\begin{document}

\title{Dimensions of tight spans}
\author{Mike Develin}
\address{Mike Develin, American Institute of Mathematics, 360 Portage 
Ave., Palo Alto, CA 94306-2244, USA}
\date{\today}
\email{develin@post.harvard.edu}

\begin{abstract}
Given a finite metric, one can construct its tight span, a geometric object representing the metric. The dimension of a tight span encodes, among other things, the size of the space of explanatory trees for that metric; for instance, if the metric is a tree metric, the dimension of the tight span is one. We show that the dimension of the tight span of a generic metric is between $\lceil \frac{n}{3}\rceil$ and $\lfloor \frac{n}{2}\rfloor$, and that both bounds are tight.
\end{abstract}

\maketitle

\section{Introduction}

Let $d$ be a metric on a set of $n$ points labeled with the elements of $[n] := \{1,\ldots,n\}$, i.e. a function
\[
d:{n\choose 2}\rightarrow \RR^+\;{\rm such}\;{\rm that}\; d_{ij}+d_{jk}\ge 
d_{ik}\;{\rm for}\;{\rm all}\;i,j,k\in [n].
\]

The {\rm injective hull}~\cite{Isb} or {\rm tight span} of the metric is a
geometric object encoding it, generalizing the corresponding tree in the 
case of a tree metric~\cite{Dress}; Dress, 
Huber, and Moulton~\cite{DHM} showed 
that it is given by
the complex of bounded faces of the polyhedron
\[
P_d = \{x\in (\RR^+)^n\,\mid\,x_i+x_j\ge d_{ij}, 1\le i<j\le n\}.
\]
In~\cite{SY}, Sturmfels and Yu observed that this is polar to the complex 
of 
interior faces of the regular subdivision $\Delta_d$ of the hypersimplex
\[
\Delta(n,2) :=\;{\rm conv}\;\{e_i+e_j\,\mid\, 1\le i<j\le n\}\subset \RR^n
\]
given by lifting a vertex $e_i+e_j$ to height $d_{ij}$ and taking the complex of upper faces of the resulting polytope. 

Another formulation of the tight span mentioned in~\cite{SY} is the following. A metric corresponds to assigning a weight $d_{ij}$ to each edge $ij$ of the complete graph $K_n$; a subgraph $G$ corresponds to a cell of $\Delta_d$ if there exists an $x\in \RR^n$ satisfying
\begin{eqnarray*}
x_i+x_j &=& d_{ij}\;{\rm if}\; ij\in G \\
x_i+x_j &>& d_{ij}\;{\rm if}\; ij\notin G.
\end{eqnarray*}
The vector $(x,1)$ gives the coordinates of the hyperplane defining the 
corresponding upper cell in the lifted hypersimplex.

As a polyhedral complex, the tight span has a dimension. This dimension measures the {\em combinatorial dimension} of the metric; for instance, if the metric is a tree metric (so that it can be realized by placing the $n$ points on a tree with weighted edges and taking the resulting distances), the dimension of the tight span will be one. In a sense, the dimension measures how far the metric is from being a tree metric.

It is easy to express dimension in terms of the hypersimplex formulation, due to the equivalence of complexes outlined above. In particular, the dimension of the tight span is equal to the maximal codimension of any interior cell of the corresponding triangulation. 

This dimension for arbitrary metrics can be any number between one and $\lfloor \frac{n}{2}\rfloor$. It makes sense to restrict to metrics satisfying the following genericity condition.

\begin{defn}[\cite{SY}]
A metric is {\em generic} if each cell of the corresponding subdivision of $\Delta(n,2)$ is a simplex.
\end{defn}

For instance, this definition forbids such non-generic behavior as 
$d_{12}+d_{34}=d_{13}+d_{24}$, which corresponds to four points alignable 
on a tree. The genericity condition consists of being in a 
full-dimensional cell of the {\em metric fan}~\cite{SY} which partitions 
all metrics into combinatorial equivalence classes of tight spans.

Computation up to $n=4$ is easily done by hand; all generic metrics on 2
or 3 points have combinatorial dimension one, while all generic metrics on
four points have combinatorial dimension two. For $n=5$, all generic
metrics have combinatorial dimension two~\cite{Dress}, but for $n=6$, 
Sturmfels and Yu
showed that, surprisingly, generic metrics do not all have the same
dimension; they can have dimension two or three.

In this paper, we complete the classification of combinatorial dimensions of metrics, proving the following theorem.

\begin{thm}\label{mainthm}
The combinatorial dimension of any generic metric on $n$ points lies between $\lceil\frac{n}{3}\rceil$ and $\lfloor\frac{n}{2}\rfloor$, inclusive. Both bounds are tight.
\end{thm}

Along the way, we present a connection to integrality of a certain linear program, and develop a theory of what the faces in a triangulation of a hypersimplex look like. Note that Theorem~\ref{mainthm} implies a corresponding result about triangulations of hypersimplices, namely that any such triangulation has its smallest interior cell of codimension between $\lceil\frac{n}{3}\rceil$ and $\lfloor\frac{n}{2}\rfloor$.

\section{Linear Programming}

Take a generic finite $n$-point metric $d$ in its graph representation, i.e. an edge-labeled $K_n$ with labels satisfying the triangle inequality. Let $0\le \omega\in \RR^n$ be a nonnegative vector (in practice, usually integral), and denote by $|\omega|$ half the $L^1$-norm of $\omega$, $\frac{\sum \omega_i}{2}$. Then we make the following definitions.

\begin{defn}
A {\em fractional $\omega$-matching} $c$ is an assignment of a weight $c_{ij}\ge 0$ to each edge of 
$K_n$ such 
that $\omega_i = \sum_j c_{ij}$ for all $i$. If $\omega = (1,\ldots, 1)$, we call this a {\em fractional 1-matching}. The {\em support} supp$(c)$ of a fractional $\omega$-matching is the set of edges $ij$ with $c_{ij}>0$.
\end{defn}

Basically, $\omega$ represents the desired valences of the vertices, and a fractional $\omega$-matching is an assignment of (possibly non-integral) multiplicities to the edges to create a graph with those valences. We will be especially concerned with the following fractional $\omega$-matchings.

\begin{defn}
A fractional $\omega$-matching $c$ is called {\em LP-optimal} if $c\cdot d = \sum c_{ij} d_{ij}$ is maximal among all fractional $\omega$-matchings.
\end{defn}

These LP-optimal fractional $\omega$-matchings are the key actors in our investigation. For one, they correspond to cells in the corresponding subdivision.

\begin{prop}
Let $c$ be an LP-optimal fractional $\omega$-matching. Then supp$(c)$ is a cell in the corresponding subdivision. Furthermore, all cells arise in this fashion.
\end{prop}

\begin{proof}
Consider the point $x = \frac{\omega}{|\omega|}\in \Delta(n,2)$. In the lifted hypersimplex, the points $(x,h)$ come from taking convex combinations of lifted vertices to get $x$; these are precisely fractional $\frac{\omega}{|\omega|}$-matchings.

Since $c$ is LP-optimal, $\frac{c}{|\omega|}$ is also, and thus it maximizes $h$ over all such linear combinations. Therefore, its support set is a subset of the face containing $x$ in its interior. Since $d$ is generic, this face is a simplex, so supp$(c)$ is itself a cell as desired.

For the converse direction, take any cell, and take $\omega$ a point in its relative interior. Then we can write $\omega$ as a convex combination of all of its vertices, and the corresponding fractional $\omega$-matching given by the coefficients of this convex combination will be LP-optimal, since the lifted relevant vertices are all in the same face of the upper envelope, meaning that combining them yields a point $(\omega, h)$ in the upper envelope. This completes the proof.
\end{proof}

Therefore, we only need to investigate LP-optimal fractional $\omega$-matchings in order to determine the cells of the corresponding subdivision. We now prove a series of propositions connecting the two. Our first two propositions connect these matchings to the problem we are using them to solve.

\begin{prop}\label{interior}
Let $G$ be the support of an LP-optimal fractional $\omega$-matching. Then the corresponding cell is 
interior if and only if $G$ is a spanning subgraph of $K_n$ which is not a $K_{1,n-1}$.
\end{prop}

\begin{proof}

If $G$ is not a spanning subgraph, then there exists some vertex $i$ which is not in any of its edges.
The corresponding cell then lives in the facet-defining hyperplane $x_i = 0$ of the hypersimplex. 
Similarly, if $G$ is a $K_{1,n-1}$ with the special element $i$, then the corresponding cell lives in 
the facet-defining hyperplane $x_i=1$ of the hypersimplex.

Conversely, if $G$ is a spanning subgraph which is not a $K_{1,n-1}$, then taking the average of the 
vertices of the corresponding
cell, we obtain a point in the cell with all entries different from zero and one, which lies in the 
interior of the
hypersimplex. Thus the cell must be an interior cell.

\end{proof}

\begin{prop}
Let $G$ be the support of an LP-optimal fractional $\omega$-matching. The dimension of the corresponding cell of the triangulation of $\Delta(n,2)$ is one less than the number of edges of $G$.
\end{prop}

\begin{proof}
The number of vertices of the corresponding cell is equal to the number of edges of $G$ by definition. Since the metric is generic, all cells are simplices, so the dimension of this cell is that number minus one.
\end{proof}

We now know what statistic to investigate to solve the problem: we need to determine the minimum number of edges in a spanning support of an LP-optimal fractional $\omega$-matching. We now embark upon the investigation of such graphs.

\begin{prop}\label{noeventours}
Let $G$ be the support of an LP-optimal fractional $\omega$-matching. Then $G$ has no nontrivial even tours, where a tour is a sequence of edges $i_k j_k$ with $j_k = i_{k+1}$ beginning and ending at the same vertex, and a tour is trivial if every edge appears at least twice.
\end{prop}

\begin{proof}
Suppose $G$ had a nontrivial even tour. Let the corresponding vertices of the hypersimplex be $v_1, w_1, \ldots, v_n, w_n$. Then $v_1 + \ldots + v_{n} = w_1 + \ldots + w_n$ is a nontrivial affine dependence among the vertices of the corresponding cell, so it this cell is not a simplex, contradicting the fact that $d$ is generic.
\end{proof}

Therefore, all connected components of such a graph must either be trees or have exactly one cycle of
odd length.

\begin{prop}\label{unique}
Given a vector $\omega$, there exists a unique LP-optimal fractional $\omega$-matching.
\end{prop}

\begin{proof}
Suppose there are two distinct LP-optimal fractional $\omega$-matchings $c_1$ and $c_2$. Then we can 
write the point $\frac{\omega}{|\omega|}$ in two different convex combinations, each of which lifts to 
an upper-envelope convex combination. This yields an affine dependence among the vertices of the face 
given by the LP-optimal fractional $\omega$-matching $\frac{c_1+c_2}{2}$ (i.e. 
supp$(c_1)\,\cup\,$supp$(c_2)$), a contradiction since it would then not be a simplex and $d$ is 
generic.
\end{proof}

This means that any $\omega$ corresponds to a unique LP-optimal fractional $\omega$-matching. The next 
proposition establishes a sort-of converse.

\begin{prop}\label{suppenough}
Let $G$ be the support of an LP-optimal fractional $\phi$-matching $r$. Then all fractional 
$\omega$-matchings (for all $\omega$) with support $G$ are LP-optimal.
\end{prop}

\begin{proof}
Suppose not. Then we have some $c$ a fractional $\omega$-matching with support $G$ and $c^\prime$ a 
fractional $\omega$-matching with $c^\prime \cdot d > c\cdot d$. Take $\epsilon$ small enough so that 
$\epsilon\cdot c < r$. Then $\epsilon c^\prime \cdot d > \epsilon c \cdot d$, so $(\epsilon c^\prime - 
\epsilon c + r) \cdot d > r \cdot d$. Since $\epsilon c^\prime$ and $\epsilon c$ are both 
$\omega$-matchings, $\epsilon c^\prime - \epsilon c + r$ is a $\phi$ matching, contradicting 
LP-optimality of $r$. 
\end{proof}

In the hypersimplex language, this is saying that if some convex combination of lifted points is in the 
upper envelope, then every convex combination of those points is, a routine statement. Our next 
proposition gives the associated linear program.

\begin{prop}
Let $G$ be any graph, and let $\omega$ be given by letting $\omega_i = \text{deg}_i (G)$. Then $G$ 
corresponds to a cell if and only if the linear program given by maximizing $c\cdot d$ on the polytope 
given by 
\[
c_{ij}\ge 0\;\forall i,j\,; \sum c_{ij} = \omega_i\;\forall i
\]
has (the indicator function of) $G$ as its optimal vertex.
\end{prop}

\begin{proof}
$G$ corresponds to a cell if and only if there is some $\phi$-matching $r$
with support $G$. By Proposition~\ref{unique}, this holds if and only if
$G$ itself is LP-optimal, which amounts to the statement of the theorem.
\end{proof}

The final proposition of this section presents a criterion for genericity in terms of LP-optimality.

\begin{prop}\label{genericcond}
A metric $d$ is generic if and only if its supports of LP-optimal fractional $\omega$-matchings (for 
all $\omega)$ all have no nontrivial even tours.
\end{prop}

\begin{proof}

Proposition~\ref{noeventours} provides one direction of the proof. For the
other direction, suppose $d$ is not generic. Take some signed affine
dependence $C$ among the vertices of a face, which we can write as $\omega
:= \sum_{i\in I} c_i v_i = \sum_{j\in J} b_j v_j$ with $I\cap J =
\emptyset$ and $\sum c_i = \sum b_j = 1$. Then $c$ and $b$ are both
LP-optimal fractional $\omega$-matchings; $c+b$ is also.  supp$(c)$ and
supp$(b)$ are disjoint spanning subgraphs of $K_{{\rm supp}(w)}$ (the
complete graph on the coordinates of $\omega$ which are nonzero) , so
their union, which is supp$(c+b)$, must contain a nontrivial even tour
(alternately take edges from supp$(c)$ and supp$(b)$ until we complete a
tour.)

\end{proof}

In the next section, we leverage this LP-theory to make statements about the dimension of the tight 
span 
of a generic metric $d$.

\section{Dimension of tight spans}

In this section, we prove Theorem~\ref{mainthm}. We start by proving the indicated bounds.

\begin{thm}
Let $d$ be a generic metric. Then the tight span of $d$ has dimension at least $\lceil 
\frac{n}{3}\rceil$ and at most $\lfloor \frac{n}{2}\rfloor$
\end{thm}

\begin{proof}

The dimension of the tight span is equal to the maximal codimension of an
interior simplex in the corresponding triangulation of $\Delta(n,2)$.
Since this polytope is $\lfloor\frac{n-1}{2}\rfloor$-neighborly, this
codimension can be at most $(n-1) - \lfloor \frac{n-1}{2}\rfloor = \lfloor
\frac{n}{2}\rfloor$.

For the lower bound, we use the LP-theory developed in the previous
section. Let $c$ be the LP-optimal fractional 1-matching, which is unique
by Proposition~\ref{unique}. We claim that the connected components of
supp$(c)$ are all odd cycles and isolated edges. Suppose that some
connected component of supp$(c)$ has a leaf $i$ connected to only one
vertex $i^\prime$. Since $1 = \sum c_{ij} = c_{ii^\prime}$, and $1 = 
\sum_j c_{i^\prime j}
= 1 + \sum_{j\neq i} c_{i^\prime j}$, $i^\prime$ must only be connected to 
$i$, and
that connected component is a single edge.

If a connected component has no leaves, it contains a cycle, which must be odd by
Proposition~\ref{noeventours}. If $G$ contains any other edges, since it has no leaves, it follows that 
it
must have a nontrivial even tour (follow edges not in the cycle until you reintersect the cycle, then
follow whichever half of the cycle gives you an even total length.) So that connected component must 
be just an odd cycle, completing the proof of the claim.

The edges in supp$(c)$ form a cell. We claim that we can find a subset of size at most $\lfloor 
\frac{2n}{3}\rfloor$ which spans $G$. Indeed, it is immediate that taking all isolated edges and 
$\lceil \frac{k}{2}\rceil$ spanning edges from each $k$-cycle does the trick, as each connected 
component of this subgraph has either one edge and two vertices or two edges and three vertices. By 
Proposition~\ref{interior}, this subset corresponds to an interior cell, and its codimension is 
$(n-1)-(\lfloor \frac{2n}{3}\rfloor-1) = \lceil \frac{n}{3}\rceil$ as desired.
\end{proof}

Whether a tight span has the maximum dimension of $\frac{n}{2}$, for $n$ even, corresponds to testing 
integrality of a linear program, namely whether or not the LP-optimal fractional 1-matching is 
integral. If this matching is integral, its support is a set of $\frac{n}{2}$ edges, which corresponds 
to an interior cell of codimension $\frac{n}{2}$. Conversely, an interior face of codimension 
$\frac{n}{2}$ corresponds to a matching, which must be LP-optimal as otherwise its edges do not form a 
face.

In particular, the upper bound of Theorem~\ref{mainthm} is easy to achieve.

\begin{prop}\label{matchex}
If $d$ is the metric given by $d_{ij} = 
2$ for $|i-j| = \lfloor \frac{n}{2}\rfloor$ and $d_{ij} = 1+\alpha_{ij}$ otherwise, where the 
$\alpha_{ij}$'s are positive numbers smaller than $1/n^2$ forming a transcendence 
basis over $\QQ$. Then 
$d$ is generic and 
the tight span of $d$ has 
dimension $\lfloor 
\frac{n}{2}\rfloor$. 
\end{prop}

\begin{proof}
We need to show first that $d$ is generic. By Proposition~\ref{genericcond}, it suffices to show that 
no LP-optimal fractional $\omega$-matching $c$ has a nontrivial even tour. Suppose one does; then let 
the odd edges of the tour comprise the (multi)set $O$, and the even edges $E$. We can assume that 
$O\cap E = \emptyset$, since if $O\cap E \neq \emptyset$, the tour breaks up into two even subtours. 
Furthermore, there is some $\epsilon$ such that $\epsilon < c_{ij}$ for all $i,j$ with the edge $ij$ in 
$O$ or in $E$. Since reducing each $c_{ij}$ with $ij\in O$ by $\epsilon$ and increasing each $c_{ij}$ 
with $i,j\in E$ by $\epsilon$ yields another fractional $\omega$-matching, and so does the reverse 
operation, we must have $\sum_{ij\in O} d_{ij} = \sum_{ij\in E} d_{ij}$.

However, by our choice of the $d_{ij}$'s, the only way such sums can be equal is if all relevant
$d_{ij}$'s are 2, but it is impossible to form a nontrivial even tour with only 2-edges, since they
form disjoint edges (and possibly one two-edge connected component.)

If $n$ is even, the LP-optimal 1-matching is given by the edges $ij$ for $j = i + \frac{n}{2}$.  This
is spanning and has support of the appropriate size. If $n=2k+1$ is odd, the LP-optimal 1-matching is
given by taking $\frac{1}{2}(1(1+k) + 1(1+2k) + (1+k)(1+2k))$ along with the edges $ij$ for $j = i +
k$, $2\le i\le k$. This again has codimension $\lfloor \frac{n}{2}\rfloor$. 
\end{proof}

Indeed, for $n$ even, any metric with an integral LP-optimal 1-matching will yield a metric of 
dimension $\lfloor\frac{n}{2}\rfloor$. The condition that the $\alpha_{ij}$'s form a transcendence 
basis is only needed to imply that no subset of them has the same sum as any other subset modulo 1, 
which is needed for genericity.

The lower bound
is trickier to achieve. We give a construction which shows it is tight for $n=3k$; 
the 
construction is easily modified to produce examples for $n=3k+1$ and $n=3k+2$. This metric has triples 
of points which are pairwise far apart, while the distances between points in different triples are all 
small.

\begin{prop}\label{triangex}
Let $n=3k$. Let $d$ be a metric on $[n]$ given as follows 
\begin{eqnarray*}
d_{ij} & = & 2\;{\rm if}\; \lceil \frac{i}{3}\rceil = \lceil\frac{j}{3}\rceil \\
d_{ij} & = & 1+\alpha_{ij}\;{\rm otherwise},
\end{eqnarray*}
where the $\alpha_{ij}$'s are positive numbers smaller than $1/n^2$ forming a 
transcendence basis over 
$\QQ$.
Then $d$ is generic, and the tight span of $d$ has dimension 
$\frac{n}{3}$.
\end{prop}

\begin{proof}
We need to show first that $d$ is generic. As in Proposition~\ref{matchex}, if $d$ were not generic, we 
could find a nontrivial even tour with $\sum_{ij\in O} d_{ij} = \sum_{ij\in E} d_{ij}$, where $O$ is 
the multiset of odd edges of the tour and $E$ is the multiset of even edges of the tour. Again by our 
choice of the $d_{ij}$'s, the only way such sums can be equal is if all relevant
$d_{ij}$'s are 2, but it is impossible to form a nontrivial even tour with only 2-edges, since they
form disjoint triangles.

Next, we need to show that no interior face has dimension greater than $\frac{n}{3}$, i.e. that all 
supports of LP-optimal fractional $\omega$-matchings, for $\omega>0$, have at least $\frac{2n}{3}$ 
edges.

Suppose we have such a support $G$. Then by Proposition~\ref{suppenough}, (the indicator function of) 
$G$ is itself LP-optimal, and a spanning subgraph since $\omega > 0$.  We enumerate the edges as 
follows. Let $a_i$ be the number of 2-triangles $\{3r+1, 3r+2, 3r+3\}$ containing exactly $i$ edges 
from $G$. Then $G$ contains $a_1 + 2a_2 + 3a_3$ edges from 2-triangles, and $n=3(a_0 + a_1+a_2+a_3)$. 
We now enumerate edges of $G$ not in 2-triangles.

Suppose $G$ contains no edges from a 2-triangle $\{1,2,3\}$. Then it
contains edges $1x_1, 2x_2$, and $3x_3$, but we must have $x_1=x_2$ as
otherwise we can replace the edges $1x_1$ and $2x_2$ by $12$ and $x_1x_2$,
to get a better matching, contradicting LP-optimality of $G$. Similarly, 
we
must have $x_1=x_3$. So $G$ contains three intertriangle edges incident
upon this 2-triangle, giving us a total of $3a_0$ intertriangle edges.

Similarly, if $G$ contains 1 edge from a 2-triangle $\{1,2,3\}$, let it be $12$. Then $G$ must contain 
$3x$ for some $x$, giving us $a_1$ more intertriangle edges. So assuming all $3a_0+a_1$ edges just 
enumerated are distinct, $G$ has at least $3a_0+a_1$ edges not in 2-triangles, for a total of 
$3a_0+2a_1+2a_2+3a_3 \ge \frac{2n}{3}$ edges. All that remains to be shown is that all enumerated 
intertriangle edges are distinct. 

Let $14$ be an intertriangle edge; we need to show that $G$ either contains two edges from $\{1,2,3\}$ 
or two edges from $\{4,5,6\}$. $G$ contains some edge $2x$; if $x\notin \{1,3,4\}$, then we can replace 
$14$ and $2x$ by $12$ and $4x$ in the matching to contradict LP-optimality of $G$. So $G$ contains an 
edge from $\{21,23,24\}$, and similarly an edge from $\{31,32,34\}$, $\{51,54,56\}$, and 
$\{61,64,65\}$.

If $G$ contains $23$, it can't contain edge $51$, or else we can 
replace 
$\{14,23,51\}$ by 
$\{12,13,45\}$ to contradict LP-optimality of $G$. Similarly, if $G$ 
contains $23$, it can't 
contain edge $56$, or else we can replace $\{23,14,56\}$ by $\frac{1}{2}\{12,13,23,45,46,56\}$ to 
contradict LP-optimality of $G$. So if $G$ contains $23$, it must contain $54$, and similarly $64$, so 
we are done. 

If $G$ contains $24$, then it can't contain $51$ (or we can replace $\{24,51\}$ by $\{12,45\}$) or 
$56$ 
(or else we replace $\{14,24,56\}$ by $\{12,45,56\}$), so it must contain $54$, and similarly $64$, and 
so we are done.

Similarly, if $G$ contains $34$, it must contain $54$ and $64$. The only remaining case is where $G$ 
contains $21$ and $31$, in which case it contains two edges from triangle 123, completing the final 
case of the proof.
\end{proof}

So both bounds in Theorem~\ref{mainthm} are in fact tight, and in fact we can get any dimension 
satisfying those bounds by titrating the triangle construction of Proposition~\ref{triangex} and the 
matching construction of Proposition~\ref{matchex}.

\end{document}